%%%%%%%%%%%%%%%%%%%%%%%%%%%%%%%%%%%%%%%%%%%%%%%%%%%%%%%%%%%%%%%%%%%%%%%%%%%
%% Jiang, Tan; Yau, Stephen S.-T.
%% 
%% Topological Invariance of Intersection Lattices of Arrangements in 
%%   $\Bbb{CP
%% 
%% Let $\scr A^*=\{l_1,l_2,\cdots,l_n\}$ be a line arrangement in 
%%   $\Bbb{CP}^2$, i.e., a collection of distinct lines in $\Bbb{CP}^2$. 
%%   Let $L(\scr A^*)$ be the set of all intersections of elements of $A^*$
%%   partially ordered by $X\leq Y\Leftrightarrow Y\subseteq X$. Let $M(\scr
%%   A^*)$ be $\Bbb{CP}^2-\bigcup\scr A^*$ where $\bigcup\scr A^*= 
%%   \bigcup\{l_i\colon\ 1\leq i\leq n\}$. The central problem of the 
%%   theory of arrangement of lines in $\Bbb{CP}^2$ is the relationship 
%%   between $M(\scr A^*)$ and $L(\scr A^*)$.
%% 
%% publ:  Bull. Amer. Math. Soc. (N.S.) 29(1993) no. 1
%% pp:    88-93
%% type:  Research Announcement        markup: amstex    file size: 17K
%% contact: u25147@uicvm.bitnet
%% 
%% copyright: American Math. Society copyright; see end of article
%% 
%% Include files necessary for this article: bull-ppt.tex
%% 
%%%%%%%%%%%%%%%%%%%%%%%%%%%%%%%%%%%%%%%%%%%%%%%%%%%%%%%%%%%%%%%%%%%%%%%%%%%
\input amstex 
\documentstyle{amsppt}
\input bull-ppt
\keyedby{bull409/kta}
\topmatter
\cvol{29}
\cvolyear{1993}
\cmonth{July}
\cyear{1993}
\cvolno{1}
\cpgs{88-93}
%\ratitle
\title Topological Invariance of Intersection Lattices\\
of Arrangements in $\Bbb{CP}^2$\endtitle
\author Tan Jiang and Stephen S.-T. Yau\endauthor
\shortauthor{Tan Jiang and S. S.-T. Yau}
\shorttitle{Intersection lattices of arrangements in
$\Bbb{CP}^2$}
\address Department of Mathematics, University of Illinois 
at Chicago,
Box 4348, M/C 249, Chicago, Illinois 60680\endaddress
\ml\nofrills{\it E-mail address}, Tan Jiang:\enspace 
u25147\@uicvm.bitnet
\mlpar
{\it E-mail address}, Stephen S.-T.\ Yau: 
u32790\@uicvm.bitnet
\endml
\date May 5, 1992 and, in revised form, December 31, 
1992\enddate
%\keywords{}
\thanks Research was partially supported by NSF grant 
DMS8822727\endthanks
\abstract Let $\scr A^*=\{l_1,l_2,\cdots,l_n\}$ be a line 
arrangement
in $\Bbb{CP}^2$, i.e., a collection of distinct lines in 
$\Bbb{CP}^2$.
Let $L(\scr A^*)$ be the set of all intersections of 
elements of
$A^*$ partially ordered by $X\leq Y\Leftrightarrow 
Y\subseteq X$. Let
$M(\scr A^*)$ be $\Bbb{CP}^2-\bigcup\scr A^*$ 
where $\bigcup\scr A^*=
\bigcup\{l_i\colon\ 1\leq i\leq n\}$. The central problem 
of the theory of
arrangement of lines in $\Bbb{CP}^2$ is the relationship 
between
$M(\scr A^*)$ and $L(\scr A^*)$.\endabstract
\endtopmatter

\document
\thm{Main Theorem} The topological type of $M(\scr A^*)$ 
determines $L
(\scr A^*)$.\ethm
As a corollary of this, we show that the algebra and 
homotopy type of
$M(\scr A^*)$ do not determine the topological type of 
$M(\scr A^*)$.

Let $\scr A=\{H_1,\dotsc,H_n\}$ be a central arrangement 
of hyperplanes
in $\Bbb C^3$. Let $M(\scr A)=\Bbb C^3-\bigcup\{H_i\colon\ 
1\leq i\leq n\}$.
There is a standard procedure in [8] or [10] for passing 
from arrangements
of hyperplanes in $\Bbb C^3$ to arrangements of lines in 
$\Bbb{CP}^2$.
In fact, $M(\scr A)=M(\scr A^*)\times\Bbb C^*$. The 
intersection lattice
$L(\scr A)$ is the set of all intersections of elements of 
$\scr A$
partially ordered by reversed inclusion. [9] shows $L(\scr 
A)$ completely
determines the cohomology ring $H^*(M(\scr A))$.

This result brings the relation between $L(\scr A^*)$ and 
$M(\scr A^*)$
into focus. An example question: Does $L(\scr A^*)$ 
determine the homotopy
type, topological type, and diffeomorphic type of $M(\scr 
A^*)$? Conversely,
do any latter invariants of $M(\scr A^*)$ determine 
$L(\scr A^*)$?

For a general class of projective arrangements in 
$\Bbb{CP}^2$, we have
shown $L(\scr A^*)$ determines the diffeomorphic type of 
$M(\scr A^*)$
[4, 5].

Falk introduced an algebraic invariant for $L(\scr A^*)$. 
For two particular
projective arrangements in $\Bbb{CP}^2$, he asked if they 
have isomorphic
Orlik-Solomon algebras [1, 2]. Rose and Terao produced 
such an isomorphism
[11, 10]. Then Falk showed the $M(\scr A^*)\roman{s}$ in 
his example have the same
homotopic type. In view of this example, one would like to 
ask whether
$L(\scr A^*)$ is determined by the topological type of 
$M(\scr A^*)$.
The purpose of this note is to announce an affirmative 
answer to the
above question.

Let us restate the main theorem more clearly:
\thm{Main Theorem} Let $\scr A^*_1$ and $\scr A^*_2$ be 
two projective
line arrangements in $\Bbb{CP}^2$. If $M(\scr A^*_1)$ is 
homeomorphic
to $M(\scr A^*_2)$, then $L(\scr A^*_1)$ is isomorphic to 
$L(\scr A^*_2)$.\ethm

In view of the results of Rose-Terao [11] and Falk [3], 
the following
statement follows immediately from the main theorem.

\thm{Corollary} There exist $\scr A^*_1$ and $\scr A^*_2$, 
two projective
line arrangements in $\Bbb{CP}^2$ such that $M(\scr 
A^*_1)$ and $M(\scr A
^*_2)$ have the same homotopy type and isomorphic 
cohomological algebra,
but not the same topological type.\ethm

Because of the main theorem, it makes the first question 
raised above more
interesting. In fact, we believe that the following 
conjecture is true.
\thm{Conjecture} For any projective line arrangement $\scr 
A^*$ in
$\Bbb{CP}^2$, the topological type of $M(\scr A^*)$ is 
determined by
$L(\scr A^*)$.\ethm

In order to prove the main theorem, we have to separate 
arrangements in
$\Bbb{CP}^2$ into two categories. An arrangement in 
$\Bbb{CP}^2$ is called
{\it exceptional} if one of its lines has at most two 
intersection points.
An arrangement in $\Bbb{CP}^2$ is called {\it 
nonexceptional} if every line
in the arrangement has at least three intersection points.

A {\it regular} neighborhood of an arrangement $\scr A^*$ 
in $\Bbb{CP}^2$
can be defined as follows: Choose a finite triangulation of
$\Bbb{CP}^2$ in which $\bigcup\scr A^*$ is a subcomplex. 
The closed star
of $\bigcup\scr A^*$ in the second barycenter subdivision 
of this triangulation
is then a regular neighborhood of $\scr A^*$.

The fundamental observation is the following lemma.
\thm{Lemma 1} Let $U_1$ and $U_2$ be regular neighborhoods 
of $\scr A^*_1$
and $\scr A^*_2$ respectively. If $M(\scr A^*_1)$ is 
homeomorphic to
$M(\scr A^*_2)$, then $\partial U_1$ is homotopic 
equivalent to
$\partial U_2$.\ethm

Let $\scr A^*$ be an arrangement in $\Bbb{CP}^2$. Suppose 
that $\scr A^*$
has $x_1,\dotsc,x_k\ (k\geq 0)$ as multiple intersection 
points
(i.e., multiplicity $t(x_i)\geq 3$). By blowing up 
$\Bbb{CP}^2$ at
$\{x_1,\dotsc,x_k\}$, we get a set $\widetilde{\scr A^*}$ 
of lines
in a blown-up surface $\widetilde{\Bbb{CP}^2}$. 
$\widetilde{\scr A^*}$ is
called an {\it associated} arrangement in 
$\widetilde{\Bbb{CP}^2}$ induced
by $\scr A^*$. Each pair of lines of $\widetilde{\scr 
A^*}$ intersects
at most one point. Let $U(\widetilde{\scr A^*})$ be a 
regular neighborhood
of $\widetilde{\scr A^*}$ and $K(\widetilde{\scr 
A^*})=\partial U
(\widetilde{\scr A^*})$. Thus $K(\widetilde{\scr A^*})$ is 
a plumbed
3-manifold which is homeomorphic to $K(\scr A^*)$, the 
boundary of a regular
neighborhood of $\scr A^*$ in $\Bbb{CP}^2$.

A class of 3-manifolds was classified by Waldhausen [12] 
in terms of {\it
graphs} and {\it reduced graph structures} of 3-manifolds. 
We call these
3-manifolds classified in [12] as {\it Waldhausen graph 
manifolds}.
\thm{Lemma 2} If $\scr A^*$ is a nonexceptional 
arrangement in $\Bbb{CP}^2$,
then $K(\widetilde{\scr A^*})$ is a Waldhausen graph 
manifold.\ethm

We define a graph $G(\widetilde{\scr A^*})$ of 
$\widetilde{\scr A^*}$
as follows. Let each vertex correspond to a line in 
$\widetilde{\scr A^*}$
with the weight of the self-intersection number of this 
line. Let each edge
correspond to the intersection point of two lines in 
$\widetilde{\scr A^*}$.

We state some definitions and results derived from [12, 
13]. Let $M$ and $N$
be compact orientable 3-manifolds. An isomorphism $\psi$ 
of $\pi_1(N)$
onto $\pi_1(M)$ is said to {\it respect the peripheral 
structure} if for each
boundary surface $F$ in $N$ there is a boundary surface 
$G$ of $M$
such that $\psi(i_*(\pi_1(F)))\subset R$ and $R$ is 
conjugate in $\pi_1(M)$
to $i_*(\pi_1(G))$ where $i_*$ denotes inclusion 
homomorphism.

\thm\nofrills{Theorem 3}\ {\rm(cf. [13, (6.5)]}. If $M$ 
and $N$ are two
Waldhausen graph manifolds and $\psi$ is an isomorphism 
from $\pi_1(N)$
onto $\pi_1(M)$ which respect the peripheral structure and 
$H_1(M)$ is
infinite, then there exists a homeomorphism from $N$ to 
$M$ which
induces $\psi$.\ethm

\thm\nofrills{Theorem 4}\ {\rm(cf. [12, (9.4)])}. Two 
Waldhausen graph
manifolds are homeomorphic if and only if the 
corresponding graphs
are equivalent.\ethm

Now suppose that $\scr A^*_1$ and $\scr A^*_2$ are two 
nonexceptional
arrangements in $\Bbb{CP}^2$ and $M(\scr A^*_1)$ is 
homeomorphic to
$M(\scr A^*_2)$. In view of Theorem 3 and Lemmas 1 and 2, 
we have that
$K(\widetilde{\scr A^*_1})$ is homeomorphic to
$K(\widetilde{\scr A^*_2})$. By Theorem 4 we conclude that 
there is an
isomorphism from $L(\widetilde{\scr A^*_1})$ to
$L(\widetilde{\scr A^*_2})$. This isomorphism also 
preserves weights
(i.e., self-intersection number). So the main theorem 
follows from
\thm{Theorem 5} Let $\scr A^*_1$ and $\scr A^*_2$ be two 
arrangements
in $\Bbb{CP}^2$. By blowing up their multiple points 
\RM(of multiplicity
$\geq 3$\RM), we obtain two associated arrangements 
$\widetilde{\scr A^*_1}$
and $\widetilde{\scr A^*_2}$ in some blown-up surfaces 
$\widetilde
{\Bbb{CP}^2}$. Then there exists an isomorphism from 
$L(\scr A^*_1)$
onto $L(\scr A^*_2)$ which preserves weights if and only 
if there is an
isomorphism from $L(\widetilde{\scr A^*_1})$ onto 
$L(\widetilde{\scr A^*_2})$.
\ethm

We next suppose that both $\scr A^*_1$ and $\scr A^*_2$ 
are exceptional.
Write
$$\align
\scr A^*_1&=\{H_0,H_1,\dotsc,H_p,\ H_{p+1},\dotsc,H_{p+
q}\},\tag 1\\
\scr A^*_2&=\{H_0,G_1,\dotsc,G_s,\ G_{s+1},\dotsc,G_{s+
t}\}\tag 2\endalign$$
where $H_0$ (respectively $G_0$) intersects with 
$H_1,\dotsc,H_p$
(respectively $G_1,\dotsc, G_s$) at one point and 
interacts with
$H_{p+1},\dotsc,H_{p+q}$ (respectively, $G_{s+1},\dotsc, 
G_{s+t})$ at
another point. If $M(\scr A^*_1)$ is homeomorphic to 
$M(\scr A^*_2)$,
then the Orlik-Solomon algebras associated to $A_1$ and 
$A_2$ are isomorphic.
It follows that $p+q=s+t$ and $pq=st$. So $L(A_1)$ is 
isomorphic to
$L(A_2)$.

Finally, we assume that $\scr A^*_1$ is exceptional, but 
$\scr A^*_2$ is not.
We need to show that $M(\scr A^*_1)$ is not homeomorphic to
$M(\scr A^*_2)$. There are four subcases to consider.
\rk{Case {\rm a}} $\scr A^*_1$ consists of at most three 
lines. We need to observe
only that the first betti number of $M(A)$ is precisely 
$|A|$. So we have
$b_1(M(A_1))\leq 3<b_1(M(A_2))$, and $M(\scr A^*_1)$ is 
not homeomorphic
to $M(\scr A^*_2)$.\endrk
\rk{Case {\rm b}} $\scr A^*_1$ is a pencil, and $|\scr 
A^*_1|\geq 4$. This follows
immediately from the following two lemmas.\endrk
\thm{Lemma 6} Let $\scr A^*$ be an arrangement in 
$\Bbb{CP}^2$. If
$\scr A^*$ is not a pencil \RM(i.e., $\bigcap\scr 
A^*=\varnothing$\RM)
and $|\scr A
^*|\geq 3$, then $b_3(M(\scr A))$, the third betti number 
of $M(\scr A)$,
is nonzero.\ethm
\thm{Lemma 7} Let $\scr A^*$ be an arrangement in 
$\Bbb{CP}^2$. If
$\scr A^*$ is a pencil \RM(i.e., $\bigcap\scr A^*$ is a 
point\RM), then
$b_3(M(\scr A))$, the third betti number of $M(\scr A^*)$, 
is zero.\ethm
\rk{Case {\rm c}} $\scr A^*_1$ consists of a pencil and a 
line in general position,
and $|A^*_1|\geq 4$ (see Figure 1).\endrk

\fighere{18pc}\caption{Figure 1}
By using Neumann's calculus of plumbing [7], one can show 
that $K(\scr A^*_1)$,
the boundary of a regular neighborhood of $\scr A^*_1$, is 
a reduced graph
manifold with reduced graph structure equal to empty set. 
It follows from
Lemma 1 and Theorems 3 and 4 that $M(\scr A^*_1)$ is not 
homeomorphic
to $M(\scr A^*_2)$.
\rk{Case {\rm d}} $\scr A^*_1=\{H_0,H_1,\dotsc,H_p,\ H_{p+
1},\dotsc,H_{p+q}\}$
where $\bigcap^p_{i=0}H_i$ and $H_0\cap(\bigcap^{p+q}_{i=p+
1}H_i)$ are
two different nonempty intersections, $p>1$ and $q>1$ (see 
Figure 2).\endrk

By blowing up the points $\bigcap^p_{i=1}H_i$ and 
$\bigcap^{p+q}
_{j=p+1}H_j$, we get the following pictures (Figure 3) 
where $E_1$ and $E_2$
are exceptional lines.

\fighere{14.5pc}\caption{Figure 2}

\fighere{18pc}\caption{Figure 3}

\noindent
Here the numbers in the parenthesis are the 
self-intersection numbers.
We blow down $H_0$ to a point and get the following 
pictures (see Figures 4
and 5).

The graph manifold $K(\scr A^*_1)$ is then a Waldhausen 
graph manifold with
the graph $G$. Notice that $G$ has only zero weights, 
while the graph
of $K(\scr A^*_2)$ has nonzero weight. Therefore, by Lemma 
1 and Theorems 3
and 4 again, we know that $M(\scr A^*_1)$ is not 
homeomorphic to
$M(\scr A^*_2)$.

\fighere{18pc}\caption{Figure 4}
\eject

\fighere{16.5pc}\caption{Figure 5}

\enddocument